\documentclass[12pt]{amsart}
\usepackage{amsmath, amssymb, amsthm, hyperref}
\usepackage[a4paper, margin=1in]{geometry}

\title[Recurrence Relations for \(k\)-Fold Nested Power Sums]{Recurrence Relations for \(k\)-Fold Nested Power Sums}
\author{Alexander R. Povolotsky}
\date{2025}

\begin{document}
\maketitle

\begin{abstract}
We consider the \(k\)-fold nested sum of integer powers, \(F(n,m,k)\), defined as repeated partial sums of the classical Faulhaber polynomials. We provide an explicit recurrence relation relating \(F(n,m,k)\) to sums of lower power \(m-1\) and higher nesting level \(k+1\). This identity is derived from a core algebraic relation on the binomial coefficients that form the kernel of the nested sum’s representation. 
We discuss the relevance to the 2010 paper by S. Butler and P. Karasik, “A Note on Nested Sums” (JIS, Vol. 13, Article 10.4.4), which studies nested sums of powers of integers that generalize Faulhaber-type sums.
We also discuss the equivalence to a related recurrence previously established in the context of hypersums of powers of integers by J.~L.~Cereceda.
\end{abstract}

Let \(\mathbb{Z}^+\) be the set of positive integers, and let integers \(m\ge0\) and \(k\ge1\). The \(k\)-fold nested sum of \(m\)-th powers, denoted \(F(n,m,k)\), is defined recursively by
\[
F(n,m,1) := \sum_{i=1}^n i^m,\qquad
F(n,m,k) := \sum_{t=1}^n F(t,m,k-1), \quad k\ge2.
\]
Alternatively, one has the closed-form representation
\[
F(n,m,k) = \sum_{r=1}^n \binom{n-r+k-1}{k-1}\, r^m.
\]

For integers \(n\ge1\), \(k\ge1\), and \(m\ge1\), the following identity holds:
\[
F(n,m,k) = n\,F(n,m-1,k) - k\,F(n-1,m-1,k+1).
\]

\begin{proof}
Starting from the representation
\[
F(n,m,k) = \sum_{r=1}^n \binom{n - r + k - 1}{k - 1}\, r^m,
\]
we substitute into the right-hand side:
\[
\text{RHS} = n \sum_{r=1}^n \binom{n - r + k - 1}{k - 1} r^{\,m-1}
- k \sum_{r=1}^{\,n-1} \binom{n - 1 - r + k}{k} r^{\,m-1}.
\]
We may extend the second sum to run \(r=1\) to \(n\), since \(\binom{k-1}{k}=0\). Both sums thus index \(r=1\) to \(n\). For each term fix \(r\); set \(a = n - r + k - 1\). Then one checks the algebraic identity
\[
n\,\binom{a}{k-1} - k\,\binom{a}{k} = r\,\binom{a}{k-1}.
\]
Hence each term in the expanded RHS becomes
\[
r\,\binom{n - r + k -1}{k-1} \, r^{\,m-1} = \binom{n - r + k -1}{k-1}\,r^{\,m},
\]
and summing over \(r=1\) to \(n\) yields exactly \(F(n,m,k)\). This completes the proof for all integers \(n,k,m\ge1\).
\end{proof}

For the boundary case \(m=0\) one obtains by the hockey-stick identity
\[
F(n,0,k) = \sum_{r=1}^n \binom{n - r + k - 1}{k - 1} = \binom{n + k - 1}{k}.
\]
The recurrence for \(m=0\) is not directly applicable, but one has
\[
F(n,0,k) = F(n-1,0,k) + F(n,0,k-1), \quad
F(n,0,1)=n,\quad F(0,0,k)=0.
\]

The 2010 paper by S. Butler and P. Karasik, “A Note on Nested Sums” (JIS, Vol. 13, Article 10.4.4), studies nested sums of powers of integers that generalize Faulhaber-type sums. Specifically, Butler \& Karasik define recursively
\[
S_k(n) = \sum_{i=1}^{n} S_{k-1}(i), \quad S_1(n) = \sum_{i=1}^n i^m,
\]
and analyze closed forms, generating functions, and polynomial properties of such nested sums. This paper defines essentially the same object, written as
\[
F(n,m,k) = \sum_{t=1}^n F(t,m,k-1),
\]
and then derives a new recurrence relation
\[
F(n,m,k) = n F(n,m-1,k) - k F(n-1,m-1,k+1),
\]
which expresses the nested power sum in terms of lower powers and higher nesting.

In J.~L.~Cereceda’s work (e.g., “Hypersums of powers of integers via the Stolz-Cesàro lemma”, Int.\ J. Pure Appl.\ Math.\ 96(3), 2014) one finds a related recurrence for hypersums of powers of integers:
\[
F(n,m,k) = \frac{k}{\,n + k\,}\,F(n,m,k+1)
+ \frac{1}{\,n + k\,}\,F(n,m+1,k), \quad n\ge1,\,k\ge1,\,m\ge0.
\]
A straightforward telescoping argument shows that this relation is equivalent to the stated theorem. Indeed, observe that the nested sums satisfy the elementary difference identity
\[
F(n,m,k) - F(n-1,m,k) = F(n,m,k-1),
\]
which expresses each level-\(k\) sum as a discrete accumulation of the level-\((k-1)\) sums. Substituting this relation into Cereceda’s recurrence and simplifying yields the form of Theorem 1. Conceptually, this “telescoping” transformation collapses adjacent nesting levels by expressing higher-order hypersums as successive differences of lower-order ones, thereby showing that both recurrences encode the same structural dependence between the indices \(m\) and \(k\).

Some instances of \(F(n,m,k)\) correspond to well-known OEIS integer sequences:
\[
\begin{aligned}
m=1,\,k=1: &\quad F(n,1,1) = \frac{n(n+1)}{2}, \\
m=1,\,k=2: &\quad F(n,1,2) = \frac{n(n+1)(n+2)}{6} \to \text{A000292}, \\
m=1,\,k=3: &\quad F(n,1,3) = \frac{n(n+1)(n+2)(n+3)}{24} \to \text{A000332}, \\
m=3,\,k=1: &\quad F(n,3,1) = \left(\frac{n(n+1)}{2}\right)^2 \to \text{A000537}.
\end{aligned}
\]


\begin{thebibliography}{9}
\bibitem{ButlerKarasik2010}
S.~Butler and P.~Karasik,
\emph{A Note on Nested Sums},
Journal of Integer Sequences, Vol. 13 (2010), Article 10.4.4.

\bibitem{Cereceda2014}
J.~L.~Cereceda,
\emph{Hypersums of powers of integers via the Stolz-Cesàro lemma},
International Journal of Pure and Applied Mathematics, 96(3), 2014, 343–351.
\end{thebibliography}
\end{document}